\input amstex
\documentstyle{amsppt}
\NoRunningHeads
\refstyle{A}
\widestnumber\no{999}
\TagsOnRight
\topmatter

\title CARTAN CALCULUS AND ITS GENERALIZATIONS \\ 
VIA A PATH-INTEGRAL APPROACH TO CLASSICAL MECHANICS \endtitle

\author Ennio Gozzi \endauthor

\address Dipartimento di Fisica Teorica, Universita' di Trieste,
 Strada Costiera 11,\linebreak
 Miramarare-Grignano 34014 Trieste 
and INFN, sezione di Trieste, Italy \endaddress
\thanks This work has been partly supported by  a grant
from the Ministero della ricerca scientifica e tecnologica of Italy  
and by a Nato grant. \endthanks
\abstract In this paper we review the recently proposed path-integral
counterpart of the Koopman-von Neumann operatorial approach to classical
Hamiltonian mechanics. We identify in particular the geometrical variables
entering this formulation and show that they are essentially a basis
of the cotangent bundle to the tangent bundle to phase-space. In this
space we introduce an extended Poisson brackets structure which allows us
to re-do all the usual Cartan calculus on symplectic manifolds via these
brackets. We also briefly sketch how the Schouten-Nijenhuis, the Fr\"olicher-
Nijenhuis and the Nijenhuis-Richardson brackets look in our formalism.

\endabstract
\endtopmatter
\document
\specialhead 1. Introduction \endspecialhead
Soon after the appearance of quantum mechanics
with its intrinsic operatorial structure, Koopman
and von Neumann \cite{6}\cite{11} gave an operatorial formulation also to classical
Hamiltonian mechanics.
\par 
It is well-known that any theory which exists in the operatorial formulation
admits also a path-integral version as it was shown for
quantum mechanics \cite{2} long ago by R.P.Feynman. Having, thanks to the
work of Koopman and von Neumann, also classical
mechanics in an operatorial version , it was easy to
give a path-integral for it \cite{4}. The weight in this classical
path-integral (CPI) is just a functional Dirac delta forcing all paths on
the classical ones.
\par
We will briefly review ref.[4] in section 2 showing how the simple
Dirac delta mentioned above can be turned into an effective evolution
classical operator. In  section 3 We will illustrate  the
geometrical meaning of the various variables appearing in our path-integral,
variables which do not parametrize only the phase space ~${\Cal M}$~
of the system but instead the ~$T^{\ast}(T{\Cal M})$~ which is the
cotangent bundle to the tangent bundle to phase-space. In this space 
(being a cotangent bundle) there naturally exists an extended  Poisson 
structure ({\it epb}). In the same section we will then show how one can 
reproduce all the operations of the
usual Cartan calculus on symplectic manifolds via our {\it epb}
and via some universal charges present in our CPI.
We conclude the paper with section 4 where the Schouten-Nijenhuis (NS)\cite{7}
\cite{10},
the Fr\"olicher-Nijenhuis (FN)\cite{3} and the Nijenhuis-
Richardson (NR)\cite{8} brackets are
built out of our {\it epb} and the variables of ~$T^{\ast}(T{\Cal M})$.

\specialhead 2. Classical Path-Integral \endspecialhead

We shall  briefly review in this section the path-integral formulation
of classical mechanics\cite{4}. 
The propagator ~$P\bigl(\phi_{2},t_{2}\vert \phi_{1},t_{1}\bigr)$, which gives
 the {\it classical}
probability for a particle to be at the point
~$\phi_{2}$~at time ~$t_{2}$, given that it was at the point~$\phi_{1}$~at
time $t_{1}$, is just a delta function
$$P\bigl(\phi_{2},t_{2}\vert
\phi_{1},t_{1}\bigr)={\delta}^{2n}\bigl(\phi_{2}-\Phi_{cl}(t_{2},
\phi_{1})\bigr)\tag{2-1}$$
where ~$\Phi_{cl}(t,\phi_{0})$~ is a solution of Hamilton's equation
${\dot\phi}^{a}(t)=\omega^{ab}\partial_{b}H(\phi(t))$
subject to the initial conditions ~$\phi^{a}(t_{1})=\phi_{1}^{a}$~
Here ~$H$~ is the conventional Hamiltonian of a dynamical system
defined on some phase-space ~${\Cal M}_{2n}$ with local coordinates\linebreak
$\phi^{a},a=1\cdots 2n$~ and a constant symplectic structure~
$\omega={1\over 2}\omega_{ab}d\phi^{a}\wedge d\phi^{b}$.
\par
The delta function in ~\thetag{2-1}~ can 
be rewritten\footnote{We will often write ~$\phi$~without putting
the upper indices ~$a$. The lower indices will instead indicate if
they are the first or last point of a trajectory.} as
$$\delta^{2n}\bigl(\phi_{2}-\Phi_{cl}(t_{2},\phi_{1})\bigr)=
\Bigl\{\prod_{i=1}^{
 N-1}\int d\phi_{(i)}
 \delta^{2n}\bigl(\phi_{(i)}-\Phi_{cl}(t_{i},\phi_{1})\bigr)\Bigr\}
\delta^{2n}\bigl(\phi_{2}-\Phi_{cl}(t_{2},\phi_{1})\bigr)
\tag{2-2}$$
where We have sliced the interval  [0,t] in N intervals and
labelled the various instants as ~$t_{i}$ and the fields at ~$t_{i}$~
as
~$\phi_{(i)}$. Each delta function
contained in the product on the RHS of ~\thetag{2-2} can be written as:
$$\delta^{2n}\bigl(\phi_{(i)}-\Phi_{cl}(t_{i},\phi_{1})\bigr)=\prod_{a=1}^{2n}
\delta\bigl({\dot\phi}^{a}-\omega^{ab}\partial_{b}H\bigr)_{\vert
t_{i}}
det\bigl[\delta^{a}_{b}\partial_{t}-\partial_{b}\bigl(\omega_{ac}
(\phi)\partial_{c}H(\phi)\bigr)\bigr]_{\vert t_{i}}\tag{2-3}$$
where the argument of the determinant is obtained from the functional
derivative of the equation of motion with respect
to~$\phi_{(i)}$. Introducing Grassmannian variables ~$c^{a}$~and~
${\bar c}_{a}$~ to exponentiate the determinant\cite{9}, and an
auxiliary variable ~$\lambda_{a}$~ to exponentiate the delta
functions, one can re-write the propagator above as a path-integral.
$$P\bigl(\phi_{2},t_{2}\vert
\phi_{1},t_{1}\bigr)=\int_{\phi_{1}}^{\phi_{2}}
{\Cal D}\phi~{\Cal D}\lambda~
{\Cal D}c~{\Cal D}{\bar c}~exp~i{\widetilde S}\tag{2-4}$$
where~
${\widetilde S}=\int^{t_{2}}_{t_{1}}dt~{\widetilde{\Cal L}}$
with
$${\widetilde{\Cal L}}\equiv
\lambda_{a}\bigl[{\dot\phi}^{a}-\omega^{ab}\partial_{b}H(\phi)\bigr]+
i{\bar c}_{a}\bigl(\delta^{a}_{b}\partial_{t}-\partial_{b}
[\omega^{ac}\partial_{c}H(\phi)]\bigr)c^{b}\tag{2-5}$$
In the path-integral~\thetag{2-4} we have used the slicing~\thetag{2-2}
~and then taken the limit  of~$N\rightarrow\infty$. 
Holding~$\phi$~ and ~$c$~ both fixed at the endpoints of the
path-integral, one obtains the kernel,~$K(\phi_{2},
c_{2},t_{2}\vert
\phi_{1},c_{1},t_{1})$, which propagates distributions in the space
~$(\phi,c)$
$${\widetilde\varrho}(\phi_{2},c_{2},t_{2})=\int d^{2n}\phi_{1}~
d^{2n}c_{1}~{K}(\phi_{2},c_{2},t_{2}\vert \phi_{1},c_{1},t_{1})
{\widetilde \varrho}(\phi_{1},c_{1},t_{1})\tag{2-6}$$
The distributions ~${\widetilde\varrho}(\phi,c)$~ are finite sums
of monomials of the type
$${\widetilde \varrho}(\phi,c)={1\over p!}{\varrho}^{(p)}_{a_{1}\cdots
a_{p}}(\phi)~c^{a_{1}}\cdots c^{a_{p}}\tag{2-7}$$
The kernel ~${K}(\cdot\vert\cdot)$ is represented by the
path-integral
$${K}(\phi_{2},c_{2},t_{2}\vert \phi_{1},c_{1},t_{1})=\int
{\Cal D}\phi^{a}{\Cal D}\lambda_{a}{\Cal D} c^{a}
{\Cal D}{\bar c}_{a}~exp~i\int_{t_{1}}^{t_{2}}dt
{\widetilde {\Cal L}}\tag{2-8}$$
with the boundary conditions ~$\phi^{a}(t_{1,2})=\phi^{a}_{1,2}$
and ~$c^{a}(t_{1,2})=c^{a}_{1,2}$. 
It is also easy from here to build a classical generating functional
~$Z_{cl}$~ from which all correlation-functions can be derived.
It is given by
$$Z_{cl}=\int {\Cal D}\phi^{a}(t)~{\Cal D}\lambda_{a}(t)~{\Cal D}
c^{a}(t)~{\Cal D}{\bar c}_{a}~exp~i\int dt \bigl\{{\widetilde{\Cal L}}+
source~terms\bigr\}\tag{2-9}$$
where the Lagrangian can be written as
~${\widetilde{\Cal L}}=\lambda_{a}{\dot{\phi}}^{a}+i{\bar c}_{a}{\dot c}
^{b}-{\widetilde{\Cal H}}$~with the "Hamiltonian" given by
$${\widetilde{\Cal H}}=\lambda_{a}h^{a}+i{\bar c}_{a}\partial_{b}h^{a}
c^{b}\tag{2-10}$$
and where ~$h^{a}$ are the components of the Hamiltonian vector field\cite{1}
\linebreak $h^{a}(\phi)\equiv\omega^{ab}\partial_{b}H(\phi)$.
From the path-integral ~\thetag{2-8}~ and ~\thetag{2-9}~ one can
easily check\cite{4} that
the variables~$(\phi,\lambda)$  and ~$(c,{\bar c})$~ form
conjugate pairs satisfying the ($Z_{2}$-graded) commutation
relations:

$$\aligned \bigl[\phi^{a},\lambda_{b}\bigr] &= i\delta^{a}_{b} \\
\bigl[c^{a},{\bar c}_{b}\bigr] &= \delta_{b}^{a}\endaligned\tag{2-11}
$$
The commutators above are the usual split-time commutators which one
can define in any path-integral\cite{2}. For more details see ref.[4].
Because of these commutators,  the variables
~$\lambda_{a}$~ and ~${\bar c}_{a}$~ can be represented, 
in a sort of "Schroedinger-like" picture, as 
~$$\lambda_{a}=-i{\partial\over\partial\phi^{a}}\equiv
-i\partial_{a}~~~;~~~{\bar c}_{a}={\partial\over\partial c^{a}}\tag{2-12}$$ 
So one sees immediately that the ~$\lambda_{a}$~represent a basis
in the tangent space ~$T_{\phi}{\Cal M}$. Inserting~\thetag{2-12}~ in 
~\thetag{2-10}~the Hamiltonian  becomes an operator\footnote{We will see in
section 3 that this is the Lie-derivative of the Hamiltonian flow\cite{1}.}

$${\widehat{\Cal H}}=-il_{h}\equiv h^{a}\partial_{a}+c^{b}(\partial_{b} h^{a})
{\partial\over \partial c^{a}}\tag{2-13}$$
The non-Grassmannian part of this operator coincides with the 
Liouvillian \linebreak ${\hat L}=h^{a}\partial_{a}$, which gives the evolution 
of standard distributions
~$\varrho^{(0)}(\phi)$~ in phase-space:
$$\partial_{t}\varrho^{(0)}(\phi,t)=-l_{h}\varrho^{(0)}(\phi,t)=
-{\hat L}\varrho^{(0)}(\phi,t)\tag{2-14}$$
This is the standard operatorial version of CM of Koopman
and von Neumann\cite{6}\linebreak\cite{11}. This proves that our 
path-integral is  really what is behind this operatorial formulation.

\specialhead 3. Cartan Calculus \endspecialhead

The reader at this point may start wondering what is the full
~${\widetilde{\Cal H}}$~ of eq.\thetag{2-10} with the Grassmannian part
included. In order to answer that question we have to understand the geometrical
meaning of the Grassmannian variables ~$c^{a}$ appearing in our path-integral
\thetag{2-9}. It is easy to see\cite{4} that under the time evolution
the variables ~$c^{a}$~ transform as
$$c^{\prime a}=c^{a}-c^{b}\partial_{b}h^{a}\Delta t\tag{3-1}$$
i.e. they transform as forms. So we can say that each ~$c^{a}(\phi)$~ belongs to the
cotangent fiber in ~$\phi$~ to ~${\Cal M}$, i.e. to ~$T^{\ast}_{\vert\phi}
{\Cal M}$. So the whole set of ~$c^{a}$~and~$\phi^{a}$~ make up the cotangent 
bundle ~$T^{\ast}{\Cal M}$~to~${\Cal M}$. Having realized that, let us now see  what the
other variables ~$\lambda_{a},{\bar c}_{a}$,~entering our path-integral, are.
Looking at eq.\thetag{2-12}, we notice that they are a basis of the tangent
space ($T$) to the~$(\phi,c)$-space which is ~$T^{\ast}{\Cal M}$.
We can then say that the 8n variables ~$(\phi^{a},c^{a},{\lambda}_{a},{\bar c}_{a})$
are the coordinates of ~$T(T^{\ast}{\Cal M})$~ which the tangent bundle
to the cotangent bundle to phase-space.
\par
As we know\cite{1} that ~$T(T^{\ast}{\Cal M})\backsim T^{\ast}(T{\Cal M})$~
and this last one is a cotangent bundle, we expect that there is an
extended Poisson structure (epb) on ~$T(T^{\ast}{\Cal M})$. There is in
 fact one
which is 
$$\{\phi^{a},\lambda_{b}\}_{epb}=\delta^{a}_{b}~~,~~\{c^{a},{\bar c}_{b}\}=
-i\delta^{a}_{b}~~,~~all~ others=0\tag{3-2}$$
Note that these are different from the normal Poisson brackets
on~${\Cal M}$~ which were\linebreak ~$\{\phi^{a},\phi^{b}\}_{pb}=\omega^{ab}$.
Via the extended Poisson brackets\thetag{3-2} we obtain from ~
${\widetilde{\Cal H}}$~ of (2-10) the same equations of motion as those which one
would obtain from the lagrangian ~${\widetilde{\Cal L}}$~of (2-5).
For ~$\phi^{a}$~ in particular we have that the same equations provided by
H via the normal Poisson brackets are also provided by ~${\widetilde{\Cal H}}$~
via the extended Poisson Brackets:
$$\{\phi^a,H\}_{pb}=\{\phi^{a},{\widetilde{\Cal H}}\}_{epb}\tag{3-3}$$
\par
Before proceeding further we should also point out that the ~
${\widetilde{\Cal H}}$~presents some universal invariance whose charges are
the following\cite{4}:

$$Q \equiv ic^{a}\lambda_{a},~~{\bar Q} \equiv i{\bar c}_{a}\omega^{ab}
\lambda_{b},~~
Q_{g} \equiv c^{a}{\bar c}_{a},~~
K \equiv {1\over 2}\omega_{ab}c^{a}c^{b},~~
{\bar K} \equiv {1\over 2}\omega^{ab}{\bar c}_{a}{\bar c}_{b}\tag{3-4}$$
The last thing to notice is that the variables~${\bar c}_{a}$~ transform
under time-evolution as:
$${\bar c}^{\prime}_{a}={\bar c}_{a}+{\bar c}_{b}\partial_{a}h^{b}\Delta t
\tag{3-5}$$
which is exactly how a basis of the vector-fields transform. Notice
that ~$\lambda_{a}$, even if it is~$-i\partial_{a}$, as indicated in eq. 
(2-12), does not transform under time evolution as a basis of the vector fields.
Its transformation is in fact:
$$\lambda^{\prime}_{a}=\bigl[\lambda_{a}+\lambda_{b}\partial_{a}h^{b}
\Delta t\bigr]+i{\bar c}_{i}\partial_{b}\partial_{a}h^{i}c^{b}\Delta t
\tag{3-6}$$
This is not in contradiction with~(2-12)~ because (3-6)  is
exactly how the derivatives ~$\partial\over\partial\phi^{a}$~
transform but when they are applied on functions of both~$\phi$~and ~$c$. More
work on this issue will appear in another paper.
\par
Having now established that ~$c^{a}$~ are  forms and
~${\bar c}_{a}$~ a basis for the vector fields, it is then natural
to build the following correspondence between forms and polinomials
in ~$c$~(which, due to the Grassmannian nature of the ~$c$, do not
need the use of the wedge product~$\wedge$), and between antisymmetric
multivectors fields and polinomials in ~${\bar c}$. We will indicate this
correspondence via a ~$\widehat{(\cdot)}$ symbol:
$$\aligned F^{(p)}={1\over p!}F_{a_{1}\cdot a_{p}}d\phi^{a_{1}}\wedge\cdots
\wedge d\phi^{a_{p}} & \Longrightarrow {\widehat F}^{(p)}={1\over p!}
F_{a_{1}\cdots  a_{p}}c^{a_{1}}\cdots\wedge c^{a_{p}}\\
v^{(p)}={1\over p!}V^{a_{1}\cdots a_{p}}\partial_{a_{1}}\wedge\cdots\wedge
\partial_{a_{p}} & \Longrightarrow {\widehat V}^{(p)}={1\over p!}
{\bar c}_{a_{1}}\cdots {\bar c}_{a_{p}}\endaligned \tag{3-7}$$
Using this correspondence it is then possible to rewrite all the
normal operations of the Cartan calculus\cite{1}, like doing an exterior
derivative on forms ~$dF$~, or doing an interior product between a vector field
and a form~$i_{v}F$~, or building the Lie-derivative of a vector field~$l_{h}$
, by just using
 polinomials in ~$c$ and ~${\bar c}$~ together with the extended Poisson
brackets structure and the charges built in (3-4). These rules, which we called
~$\{\cdot,\cdot\}_{epb}$-rules, are summarized  below :
$$\aligned dF^{(p)} & \Longrightarrow i\bigl\{Q,{\widehat F}\bigr\}_{epb}\\
i_{v}F^{(p)} & \Longrightarrow i\bigl\{{\widehat v}, {\widehat F^{(p)}}
\bigr\}_{epb}\\
l_{h}F=di_{h}F+i_{h}dF & \Longrightarrow -\bigl\{{\widetilde{\Cal H}},
{\widehat F}\bigr\}_{epb}\\				
pF^{(p)} & \Longrightarrow i\{Q_{g},{\widehat F}\}_{epb}\\
\omega(v,\cdot)\equiv v^{\flat} & \Longrightarrow i\{{\bar K},{\widehat V}
\}_{epb}\\
(df)^{\sharp} & \Longrightarrow i\{{\bar Q},f\}_{epb}\endaligned\tag{3-8}$$
where the last three operations indicated in (3-8) above are, respectively,
multiplying a form ~$F^{(p)}$~ by its degree ~$p$, mapping a vector field
~$V$~ into its associated one form~$V^{\flat}$~ via the symplectic form,
and building the associated Hamiltonian vector field ~$(df)^{\sharp}$~out
of a function ~$f$. One sees from above that the various abstract derivations of the
Cartan calculus are all implemented by  some charges acting via the 
epb-brackets. From the third relation in (3-8) one also can notice that the 
Lie-derivative of the Hamiltonian vector field of time evolution
becomes  nothing else than the~${\widetilde{\Cal H}}$~ of (2-10), thus
confirming that the weight-function of our classical path-integral,
generated by just a simple Dirac delta, is  the right {\it geometrical} object
associated to the time-evolution.
\par
The last question which our reader may be interested in getting an answer
is what becomes of the Lie-bracket\cite{1} of two vector fields ~$V$, $W$. The answer
is the following:
$$\bigl[V,W\bigr]_{lie}\Longrightarrow\varpropto \{{\widetilde{\Cal H}}_{V},
{\widehat W}\}_{epb}\varpropto\{{\widetilde{\Cal H}}_{W},{\widehat V}\}_{epb}
\tag{3-9}$$
where ~${\widehat W}=W^{a}{\bar c}_{a}$~and \
${\widehat V}=V^{a}{\bar c}_{a}$~
while ~${\widetilde{\Cal H}}_{V}=\lambda_{a}V^{a}+i{\bar c}_{a}\partial_{b}
V^{a}c^{b}$~ and\linebreak ${\widetilde{\Cal H}}_{W}=\lambda_{a}W^{a}+i{\bar c}_{a}
\partial_{b}W^{a}c^{b}$~are the analog of the Lie-derivatives associated 
respectively to the vector field ~$V$~and ~$W$.

\specialhead 4.Generalized Cartan Calculus. \endspecialhead

What we called "{\it Generalized Cartan Calculus}" is basically the following
set of brackets: the Schouten-Nijenhuis ones(NS)\cite{10}\cite{7}
between antisymmetric multivector fields, the Fr\"olicher-Nijenhuis
brackets(FN)\cite{ 8} and the  Nijenhuis-Richardson ones(NR)\cite{8} 
among vector-valued forms.
\subhead 4.A Schouten-Nijenhuis Brackets \endsubhead
These brackets are a generalization on multivector
fields of the Lie-brackets between vector fields. Following ref. [5]
and given two multivector fields ~$P\equiv X_{(1)}\wedge\cdots\wedge X_{(p)}
$~and ~$Q\equiv Y_{(1)}\wedge\cdots \wedge Y_{(q)}$~of rank respectively 
p and q, the NS-brackets among them is a multivector of rank
(p+q-1) given by:
$$\bigl[P, Q \bigr]_{(NS)}\equiv (-1)^{pq}\sum_{J=1}^{q}(-1)^{J+1}Y_{(1)}\wedge
\cdots\wedge{\hat{\hat Y}}_{(j)}\cdots\wedge Y_{q}\wedge[Y_{(j)}, P]
\tag{4.A-1}$$
where the ~${\hat{\hat Y}_{J}}$~ means that that vector-field has
been taken away, and ~$[Y_{j},P]=l_{Y_{j}}P$~is the Lie-derivative
of the vector field~$Y_{j}$~applied to the multivector ~$P$.
\par
The NS-brackets can easily be translated into our (epb)-formalism
via the rules established in section 3. The details of the
calculations will be presented elsewhere, but the final result is the
following:

$$\bigl[P,Q\bigr]_{(NS)}\Longrightarrow \propto\bigl\{L_{Y_{(1)}\cdots Y_{(q)}},
X^{a}_{(1)}{\bar c}_{a}\cdots X^{l}_{(p)}{\bar c}_{l}\bigr\}_{epb}\tag{4.A-2}$$

where ~$L_{Y_{(1)}\cdots Y_{(q)}}$~ is defined as:

$$L_{Y_{(1)}\wedge\cdots\wedge Y_{(q)}}\equiv \sum_{j=1}^{q}(-1)^{j-1}
Y^{a}_{(1)}{\bar c}_{a}\cdots {\widehat{\widehat {Y^{b}_{(j)}{\bar c}_{b}}}}\cdots
Y^{l}_{(q)}{\bar c}_{l}{\widetilde {\Cal H}_{Y_{j}}}\tag{4.A-3}$$

with~${\widetilde{\Cal H}_{Y_{(j)}}}=\lambda_{a}Y^{a}_{(j)}+i{\bar c}_{a}
\partial_{b}Y^{a}_{(j)}c^{b}$

The ~$L_{Y_{(1)}\cdots Y_{(q)}}$ above is  a generalization of the 
standard Lie-derivative.

\subhead 4.B Fr\"olicher-Nijenhuis Brackets \endsubhead
This is a bracket which associates\footnote{We use the notation of ref.5.} 
to two {\it vector-valued forms},
~$K\in\Omega^{k+1}({\Cal M},T{\Cal M})$~and~$V\in\Omega^{l+1}({\Cal M},
T{\Cal M})$~of rank respectively~$(k+1)$~ and $(l+1)$, a ($k+l+2)$~vector-valued
form
$$\bigl[K,V\bigr]_{(FN)}\in\Omega^{k+l+2}({\Cal M}, T{\Cal M})$$
Before proceeding we need to introduce some new notation\cite{5}.
First we have to generalize the notion\cite{1} of interior contraction~$i_{v}\Theta$~
of a form~$\Theta$~with a vector field ~$v$. The generalization is 
the contraction of an l-form~$\Theta$~with a vector-valued (k+1)-form ~$K$,
the result will be a (k+l)-form which can be contracted with (k+l)-vectors
~$X_{1}\cdots X_{k+l}$. Its precise definition is:
$$\multline 
\bigl(i_{K}\Theta\bigr)\bigl(X_{1},\cdots,X_{k+l}\bigr)
\, \equiv \\
\equiv {1\over (k+1)!(l-1)!}\sum_{\sigma\in S_{k+l}}sign~\sigma\Theta\biggl[
K(X_{\sigma_{1}}\cdots X_{\sigma_{k+1}}),X_{\sigma_{k+2}},\cdots,X_{\sigma_{k+l}}
\, \biggr]\endmultline\tag{4.B-1}$$
where~$\sigma$~is the set of permutation~$S_{k+l}$~of the
(k+l) vector fields~$X_{1}\cdots X_{k+l}$.
Having now the generalized interior contraction defined above, we can then
define a generalized Lie-derivative with respect to a vector-valued (k+1)-form
~$K$:
$$L_{K}\equiv i_{K}d+di_{K}\tag{4.B-2}$$
Using now (4.B-1) and (4.B-2), the FN brackets are defined\cite{5} in the
following implicit way:
$$\bigl[L_{K},L_{V}\bigr]\Theta\equiv L_{[K,V]_{(FN)}}\Theta\tag{4.B-3}$$
where~$\bigl[L_{K},L_{V}\bigr]$~ is the usual commutators among Lie-derivative
and ~$\Theta$~ is a form on which they act.
\par
Let us now find out how the FN-brackets appear in our epb-formalism.
The vector-valued forms ~$K$~and~$V$~ become
$$\aligned K & \Longrightarrow {\widehat K}\equiv K^{i}_{a\cdots k+1}
{\bar c}_{i}[c^{a}\cdots c^{k+1}]\\
V & \Longrightarrow {\widehat V}\equiv V^{j}_{a\cdots l+1}{\bar c}_{j}
[c^{a}\cdots c^{l+1}]\endaligned\tag{4.B-4}$$
Using this notation and the formulas of section 3, it is not
difficult to prove that
$$\bigl[K, V\bigr]_{(FN)}\Longrightarrow \propto\biggl\{{\widehat K},
\bigl\{{\widehat V}
,Q\bigr\}_{epb}\biggr\}_{epb}\tag{4.B-5}$$
Also the details of the above calculations will be presented in a forthcoming 
paper.

\subhead 4.C Nijenhuis-Richardson brackets \endsubhead
These are brackets also defined, as the FN ones, among (k+1) and (l+1) 
vector-valued
forms ~$K\in\Omega^{k+1}({\Cal M};T{\Cal M})$\linebreak $V\in\Omega^{l+1}({\Cal M};
T{\Cal M})$~ but whose result is a (k+l+1) vector valued-form.
Their exact definition\cite{5} is
$$\bigl[K,V\bigr]_{(NR)}\equiv i_{K}V-(-1)^{kl}i_{V}K\tag{4.C-1}$$
Here the ~$i_{K}$~and ~$i_{V}$~ are the generalized interior
contraction defined in (4.B-1).
\par
In the language of the epb-brackets the NR-brackets have 
a simple expression:
$$\bigl[K,V\bigr]_{NR}\Longrightarrow \propto\{{\widehat K},{\widehat V}
\}_{epb}\tag{4.C-2}$$
where the ~${\widehat K}$~and ~${\widehat V}$~are given in (4.B-4).  
The calculational details of this derivations will be presented elsewhere.

\specialhead 5.Conclusions \endspecialhead
The reader may wonder of what is the need of the dictionary we have created
between Cartan (and generalized) calculus and our epb-formalism. The
answer is in the fact that with our formalism we do not have to
take care of all the various numerical factors and signs and permutations 
(see 4.B-1) which one has to remember by heart in doing the standard abstract
Cartan calculus. In our case everything is taken care automatically by the
Grassmannian natures of the ~$c^{a}$~ and ~${\bar c}_{a}$~and the
graded structure of the epb-brackets. These, together with the five charges
(3-4), seem to be the central and only ingredients needed
to build all these operations. This reduction to these simple ingredients
seemed to me a thing to bring to the attention of the mathematics and physics
community in order to stimulate further investigations.
\head Acknowledgments: \endhead
I wish to thank G.Marmo for suggesting the study of the NS, FN, NR brackets
and for guidance and help throughout this work.  
A crucial discussion in january 1996 with G.Landi helped me in clearing
up my ideas on the space ~$T^{\ast}(T{\Cal M})$. Last, but not least, a special
thank to M.Reuter and M.Regini  for help on some technical points.

\enddocument